\documentclass[11pt]{amsart} 
\usepackage{latexsym}
\usepackage{amsfonts}
\usepackage{amsmath,amsthm,amssymb}

\newcommand{\cal}{\mathcal}

\def\epsilon{\varepsilon}
\def\phi{\varphi}

\def\hat{\widehat}

\newcommand{\Out}{\mbox{Out}}
\newcommand{\Aut}{\mbox{Aut}}

\newcommand{\stab}{\mbox{stab}}

\newcommand{\FN}{F_N}   

\newcommand{\CVN}{\mbox{CV}_N}
\newcommand{\CVNbar}{\overline{\mbox{CV}}_N}

\newcommand{\FM}{F_M}   

\newcommand{\R}{\mathbb R}
\newcommand{\Z}{\mathbb Z}

\def\strutdepth{\dp\strutbox}
\def \ss{\strut\vadjust{\kern-\strutdepth \sss}}
\def \sss{\vtop to \strutdepth{
\baselineskip\strutdepth\vss\llap{$\diamondsuit\;\;$}\null}}

\def\strutdepth{\dp\strutbox}
\def \sst{\strut\vadjust{\kern-\strutdepth \ssss}}
\def \ssss{\vtop to \strutdepth{
\baselineskip\strutdepth\vss\llap{$\spadesuit\;\;$}\null}}

\def\strutdepth{\dp\strutbox}
\def \ssh{\strut\vadjust{\kern-\strutdepth \sssh}}
\def \sssh{\vtop to \strutdepth{
\baselineskip\strutdepth\vss\llap{$\heartsuit\;\;$}\null}}

\def\strutdepth{\dp\strutbox}
\def \ss{\strut\vadjust{\kern-\strutdepth \sss}}
\def \sss{\vtop to \strutdepth{
\baselineskip\strutdepth\vss\llap{$\diamondsuit\;\;$}\null}}

\def\strutdepth{\dp\strutbox}
\def \sst{\strut\vadjust{\kern-\strutdepth \ssss}}
\def \ssss{\vtop to \strutdepth{
\baselineskip\strutdepth\vss\llap{$\spadesuit\;\;$}\null}}

\newtheorem{theorem}{Theorem}

\newtheorem{lemma}[theorem]{Lemma}
\newtheorem{proposition}[theorem]{Proposition}

\newtheorem{facts}[theorem]{Facts}
\theoremstyle{definition}
\newtheorem{definition}[theorem]{Definition}

\theoremstyle{remark}

\numberwithin{equation}{section}

\begin{document}

\author[M.~Lustig]{Martin Lustig}\address{LATP,
Centre de Math\'ematiques et Informatique, 
Aix-Marseille Universit\'e, 
39, rue F.~Joliot Curie, 
13453 Marseille 13, 
France}
  \email{\tt Martin.Lustig@univ-amu.fr}

\title[Tree-irreducible automorphisms of $\FN$] {Tree-irreducible automorphisms of free groups}

 
 
\maketitle 
 
 {\it To appear in: Extended Conference Abstracts, Spring 2013, CRM Documents, Centre de Recerca Matem\`atica, Bellaterra (Barcelona)}
 
 \bigskip
 
 Throughout this extended abstract of \cite{Lu} we denote by $\FN$ the non-abelian free group of finite rank $N \geq 2$. We assume familiarity of the reader with the basic terminology around automorphisms of free groups and Outer space; 
details as well as some back ground can be found for example in \cite{Lu07} or \cite{Vog}.  

\smallskip

An outer automorphism $\phi \in \Out(\FN)$ is called {\em iwip} (or {\em fully irreducible}) if no positive power $\phi^t$ fixes the conjugacy class of any non-trivial proper free factor $\FM \subset \FN$. 
Furthermore, $\phi$ is called {\em atoroidal} if no positive power $\phi^t$ fixes a non-trivial conjugacy class of elements of $\FN$. An automorphism of $\FN$ is atoroidal if and only if it is {\em hyperbolic} (i.e. the associated mapping torus group $\FN \rtimes_\phi \Z$ is Gromov-hyperbolic). An iwip $\phi$ is not atoroidal if and only if $\phi$ is {\em surface} (i.e. $\phi$ is induced by a homeomorphism of a surface with boundary). 

\smallskip

It follows directly from these definitions that $\phi \in \Out(\FN)$ is iwip if and only if $\phi^{-1}$ is iwip as well, and that 
$\phi$ is atoroidal if and only if $\phi^{-1}$ is also atoroidal.

\smallskip

Atoroidal iwip automorphisms have been widely recognized as the closest possible analogues in $\Out(\FN)$ of pseudo-Anosov mapping classes of surfaces.  In particular, the following well known properties correspond directly to known properties of pseudo-Anosov mapping classes of a closed surface $\Sigma$, if one ``translates'' $\FN$ into $\pi_1 \Sigma$, ``basis'' into ``standard generating system'', ``Outer space'' into ``Teichm\"uller space'', etc.:

\begin{facts}[\cite{BH, CHR, KL-7steps, LL4}]
\label{iwip-facts}
Let $\phi \in \Out(\FN)$ be an iwip automorphism. Then the following holds:
\begin{enumerate}
\item
There exists a {\em stretching factor} $\lambda_\phi > 1$ such that every 
conjugacy class 
$[w] \subset \FN$ 
which is not $\phi$-periodic
has uniform exponential growth 
$$|\hat{\phi^t([w])} |_{\cal A} \sim \lambda_\phi^t  C^{\cal A}_{[w]} \, .$$  
Here $\cal A$ is any ``fixed'' basis of $\FN$, $C^{\cal A}_{[w]} > 0$ is a constant which depends only on $[w]$ (and the choice of $\cal A$), $|\cdot|_{\cal A}$ denotes word length in $\cal A \cup \cal A^{-1}$, and by $\hat v$ we mean the cyclically reduced word in $\cal A$ which is given (up to cyclic permutation) by the conjugacy class $[v] \in \FN$. 

[Moreover, if $\phi$ is atoroidal then the only $\phi$-periodic conjugacy class is the trivial one.]

\smallskip
\item
$\phi$ is represented by an (absolute) train track map $f: \Gamma \to\Gamma$, with respect to some marking isomorphism $\pi_1\Gamma \cong\FN$. The non-negative transition matrix $M(f)$ (defined by the action of $f$ on the edges of $\Gamma$) is primitive, and its Perron-Frobenius eigenvalue is equal to $\lambda_\phi$.

\smallskip
\item
$\phi$ does not split over any very small graph-of-groups decomposition of $\FN$.
Equivalently, there is no very small simplicial tree $T$ with $\FN$-action by homeomorphisms such that $\phi$ ``commutes'' with some homeomorphisms $H: T \to T$. (By this we mean that there exists a representative $\Phi \in \Aut(\FN)$ of $\phi$ which satisfies $\Phi(w) H = H w: T \to T$ for every $w \in \FN$.)

\smallskip
\item
There is an (up to rescaling) unique $\R$-tree $T_+ = T_+(\phi)$, equipped with a very small action of $\FN$ by isometries, which is projectively $\phi$-invariant with stretching factor $\lambda_+ > 1$:
$$|| \phi([w])||_{T_+} = \lambda_+ ||[w]||_{T_+}$$
for any $w \in \FN$, where $||\cdot ||_{T_+}$ denotes the translation length on $T_+$.  Furthermore one has $\lambda_+ = \lambda_\phi$ (as given in (1) above).

\smallskip
\item
The tree $T_+(\phi)$ is indecomposable.  This means that for any two non-degenerate segments $I, J \subset T$ the segment $J$ is contained in the union of finitely many translates $w_i \overset{\circ}{I}$ of the interior $\overset{\circ}{I}$ of $I$, with $w_i \overset{\circ}{I} \cap w_{i+1} \overset{\circ}{I} \neq \emptyset$.

[Moreover, if $\phi$ is atoroidal, then the action of $\FN$ on $T_+(\phi)$ is free.]

\smallskip
\item
$\phi$ induces on the Thurston compactification $\CVNbar$ of Culler-Vogtmann's Outer space $\CVN$ a locally uniform North-South dynamics, where the sink is the projective class $[T_+(\phi)]$, and the source is the projective class $[T_-(\phi)]$, with $T_-(\phi) = T_+(\phi^{-1})$.

\end{enumerate}

\end{facts}

Since, as pointed out above, the assumption ``$\phi$ (atoroidal) iwip'' is equivalent to ``$\phi^{-1}$ (atoroidal) iwip'', all of the above facts are true for $\phi^{-1}$ as well.  The reader should however be warned that, contrary to the situation of pseudo-Anosov's on surfaces, the stretching factor $\lambda_- := \lambda_{\phi^{-1}}$ may well be different from $\lambda_+ = \lambda_\phi$.
For more information about the fine-structure of iwips and the subtle differences to pseudo-Anosov homeomorphisms see \cite{CH}.

\medskip

If one considers a surface $\Sigma$ with precisely one boundary component, and a pseudo-Anosov homeomorphism $h: \Sigma \to \Sigma$, then the induced outer automorphism $\phi = h_*$ on $\FN = \pi_1 \Sigma$ is iwip, but not atoroidal.  Indeed, all iwips $\phi$ which are not atoroidal arise precisely in this way.

\begin{facts}
\label{pseudo-Anosovs}
Let $\phi \in \Out(\FN)$ be induced by a pseudo-Anosov homeomorphism $h: \Sigma \to \Sigma$ of a surface $\Sigma$ with $1$ or more boundary components, via some isomorphism $\FN \cong \pi_1 \Sigma$.  Then one has:
\begin{enumerate}
\item
$\phi$ satisfies all of the properties (1) - 
(6)
 of Facts \ref{iwip-facts}.
 
 \smallskip
\item
However, if $\Sigma$ has at least 2 boundary components, then $\phi$ is not fully irreducible (= iwip).

\end{enumerate}
\end{facts}

The proof of fact (2) is obvious, since each boundary curve $C_i$ of $\Sigma$ generates a cyclic subgroup which, in the presence of at least one more boundary component $C_j \neq C_i$, is a free factor of $\pi_1 \Sigma$, and some power of $h$ maps $C_i$ to itself. 
The reader should note here that, contrary to the situation for pseudo-Anosov homeomorphisms on a surface, the negation of ``iwip'' is for a free group automorphisms  $\phi$ in general a much weaker property than the existence of a $\phi$-invariant small graph-of-groups decomposition of $\FN$.

\medskip

In view of these facts it seems clear that for many purposes the notion of iwip automorphisms is too restrictive. 
We propose here the following definition as a useful class of free group automorphisms:

\begin{definition}
\label{tree-irreducible}
An automorphism $\psi \in \Out(\FN)$ is called {\em tree-irreducible} if the following properties are satisfied:
\begin{enumerate}
\item
There is no very small graph-of-groups decomposition of $\FN$ which is preserved by $\psi$.  

\smallskip
\item
There is  an (up to rescaling) unique $\R$-tree $T_+ = T_+(\psi)$ with very small action of $\FN$ by isometries that is projectively $\psi$-invariant with stretching factor $\lambda_+ > 1$.

\smallskip
\item
The tree $T_+(\psi)$ is indecomposable.  

[However, the action of $\FN$ on $T_+(\psi)$ is in general not free.]

\smallskip
\item
The outer automorphism $\psi$ fixes every conjugacy class of any elliptic subgroup of $T$. 
Equivalently, if $x \in T$ is a point with non-trivial stabilizer $\stab_{\FN}(x) \subset \FN$, then $\phi$ fixes the $\FN$-orbit of the subgroup $\stab_{\FN}(x)$, and it induces on it the trivial outer automorphism.

\end{enumerate}
\end{definition}

It turns out
that this definition is not minimal: one can deduce property (1) as well as the uniqueness in property (2) from the remaining properties in Definition \ref{tree-irreducible}. On the other hand, the following 
is not totally obvious:

\begin{lemma}
\label{inverses}
The inverse of any tree-irreducible $\psi \in \Out(\FN)$ is also tree-irreducible.
\end{lemma}

The properties (1) and (2) (the latter slightly adjusted) of Facts \ref{iwip-facts} also hold for tree-irreducible automorphisms.
The most interesting property however, given its wide use for iwips, seems to be property (6).  It turns out that property (3) of Definition \ref{tree-irreducible} is the main ingredient in the proof given in \cite{LL4} for iwips, so that its main part can be applied
 word-by-word to tree-irreducible automorphisms as well. We obtain:

\begin{theorem}
\label{north-south}
Every tree-irreducible automorphism $\psi \in \Out(\FN)$ induces on the compactified Outer space $\CVNbar$ a locally uniform North-South dynamics.
\end{theorem}

As in the iwip case, the sink of this North-South dynamics is the projective class $[T_+(\psi)]$, and the source is the projective class $[T_-(\psi)]$, with $T_-(\psi) = T_+(\psi^{-1})$.

\medskip

The second main reason why one may consider seriously tree-irreducible automorphisms of $\FN$ comes from train-track technology:  Bestvina-Handel, in their celebrated paper \cite{BH}, have shown that every automorphism $\phi \in \Out(\FN)$ possesses a relative train-track representative, i.e. there is a graph $\Gamma$ with marking isomorphisms $\FN \cong \pi_1 \Gamma$ such that $\phi$ is induced by a self-map $f: \Gamma \to \Gamma$, and there exists an $f$-invariant filtration into not-necessarily connected subgraphs $\Gamma = \Gamma_r \supset \Gamma_{r-1} \supset \ldots \supset \Gamma_1 \supset \Gamma_0$, such that on each {\em stratum} $\Gamma_i $, considered modulo $\Gamma_{i-1}$, the restriction of $f$ is either zero, periodic, or train track with primitive transition matrix. 
After rising $f$ to a suitable positive power $f^t$ it follows 
that $\phi^t$ can be written as {\em semi-commuting product} 
(see \cite{Lu}) $\phi_s \circ \phi_{s-1} \circ \ldots \phi_1$, where each $\phi_i$ is either a Dehn twist automorphism, or else it is, up to an enlargement of the elliptic subgroups, a tree-irreducible automorphism.

\smallskip

The ``enlargement technique'' in the previous paragraph, however, is more delicate than one may suspect at first sight.  In particular, we use this technique to 
exhibit the first explicit examples of the following kind:

\begin{proposition}
\label{not-indecomposable}
There exist $\R$-trees $T$ with very small action of $\FN$ by isometries which are not indecomposable, but do also not decompose as dusted action or as very small graph-of-actions.
\end{proposition}

The interest in this subtle phenomenon comes from the fact that the opposite statement, but with the words ``very small'' before ``graph-of-actions'' erased, has been 
recently 
proved by Guirardel-Levitt  \cite{GL}.

\smallskip

Finally, in the last part of 
\cite{Lu} a technique is presented for producing 
tree-irreducible automorphisms from fully irreducible ones by ``puncturing the singularities of the expanding lamination'', in analogy to what is often done for pseudo-Anosov homeomorphisms of surfaces.

\end{document}